\documentclass{article}

\usepackage{arxiv}

\usepackage[utf8]{inputenc} % allow utf-8 input
\usepackage[T1]{fontenc}    % use 8-bit T1 fonts
\usepackage{hyperref}       % hyperlinks
\usepackage{url}            % simple URL typesetting
\usepackage{booktabs}       % professional-quality tables
\usepackage{amsfonts}       % blackboard math symbols
\usepackage{amsmath}
\usepackage{amssymb}
\usepackage{nicefrac}       % compact symbols for 1/2, etc.
\usepackage{microtype}      % microtypography
\usepackage{graphicx}
\usepackage{multirow}
\usepackage{longtable}      % multi-page per-instance result tables
\usepackage{algorithm}      % algorithm floats (labeling loop, dominance)
\usepackage{algpseudocode}  % algorithmicx pseudocode
\usepackage{xcolor}
\usepackage{listings}       % code listings (Resource concept)
\graphicspath{ {./} }

\definecolor{lstkw}{rgb}{0.10,0.10,0.55}
\definecolor{lstcmt}{rgb}{0.30,0.45,0.30}
\definecolor{lststr}{rgb}{0.55,0.10,0.10}
\lstdefinestyle{cpp}{
  language=C++,
  basicstyle=\ttfamily\footnotesize,
  keywordstyle=\color{lstkw}\bfseries,
  commentstyle=\color{lstcmt}\itshape,
  stringstyle=\color{lststr},
  morekeywords={concept,requires,constexpr,consteval,co_await,nullptr,
    std,size_t,uint32_t,uint64_t,override,noexcept},
  showstringspaces=false,
  columns=fullflexible,
  keepspaces=true,
  breaklines=true,
  frame=single,
  framesep=4pt,
  rulecolor=\color{black!30},
  aboveskip=1.0em,
  belowskip=0.6em,
  xleftmargin=2pt,
}
\title{\texorpdfstring{\texttt{bucket-graph-spprc}}{bucket-graph-spprc}: an extensible C++ library for the\\ shortest path problem with resource constraints}

\author{
 Simon Spoorendonk \\
  Denmark \\
  \texttt{simon@spoorendonk.dk} \\
  \href{https://orcid.org/0009-0007-4304-6956}{ORCID:~0009-0007-4304-6956} \\
}

\begin{document}
\maketitle

\begin{abstract}
We present \texttt{bucket-graph-spprc} (\texttt{bgspprc} for short), an open-source, header-only
C++23 library for the
shortest path problem with resource constraints (SPPRC), the pricing subproblem at the
heart of branch-cut-and-price for vehicle routing and related problems. The library
implements the bucket-graph labelling algorithm of Sadykov, Uchoa and
Pessoa~\cite{sadykov2021bucket}, with bidirectional labelling, across-arc concatenation,
bucket fixing and arc elimination, and a structure-of-arrays label store with SIMD-accelerated
dominance. Its central design feature is a compile-time \emph{resource concept}: a new SPPRC
variant is added by implementing a fixed seven-function interface, and resources compose into
a label state with no runtime dispatch, the state layout fixed at compile time. Five resources
ship built in: time/capacity, ng-path elementarity relaxation, rank-1 cuts, cumulative cost,
and pickup-and-delivery. In a reproducible, head-to-head comparison on shared public instances at
an identical bound, \texttt{bgspprc} outperforms PathWyse~\cite{salani2024pathwyse}, the main
open-source comparator, by $1.3\times$--$2.35\times$ in shifted geometric mean (and by
$1.3\times$--$2.3\times$ even when itself run single-threaded), and runs within
$1.9\times$--$2.4\times$ of parallel pull labelling~\cite{petersen2025pulllabelling}, a different
labelling technique for the same problem. The library, benchmark scripts, and pinned instances are publicly available.

\end{abstract}


\section{Introduction}
\label{sec:intro}

The shortest path problem with resource constraints (SPPRC) is the workhorse pricing
subproblem in column-generation and branch-cut-and-price approaches to vehicle routing and
related problems~\cite{pessoa2020vrpsolver}. Given a directed graph with arc costs and one or
more resources that accumulate along a path subject to per-vertex windows, the SPPRC asks for
a minimum-cost source--sink path that respects every resource constraint. In a column-generation
context the arc costs are reduced costs, so the subproblem must be solved to optimality, many
times, and is typically the computational bottleneck of the overall method. Decades of work on
labelling algorithms (dominance-based dynamic programming~\cite{desrochers1992vrptw,irnich2005survey},
elementarity relaxations~\cite{feillet2004espp,baldacci2011ng}, bidirectional
search~\cite{righini2006bidir}, and acceleration by completion bounds~\cite{sadykov2021bucket}) have
made the SPPRC tractable at the scales modern routing solvers require.

The most effective recent development is the \emph{bucket-graph} labelling algorithm of Sadykov,
Uchoa and Pessoa~\cite{sadykov2021bucket}, which partitions labels into buckets along the main
resources and uses completion bounds to fix buckets and eliminate arcs, drastically pruning the
search. It is a core ingredient of the state-of-the-art VRPSolver~\cite{pessoa2020vrpsolver}.
A subsequent line of work~\cite{sadykov2026metasolver} casts the family of SPPRC variants
encountered in routing (elementarity, rank-1 cuts, cumulative cost, pickup-and-delivery)
behind a small \emph{resource interface}, so that one labelling engine serves many problems.

Despite the maturity of these methods, the open-source landscape is thin. VRPSolver is
academic and closed. PathWyse~\cite{salani2024pathwyse} is an open (GPLv3) and well-engineered
library, but it implements \emph{standard} labelling with decremental state-space relaxation
(DSSR)~\cite{righini2008dssr}, not the bucket-graph algorithm, and resources are added at runtime by subclassing a
resource class. cspy~\cite{cspy} is a Python package aimed at flexibility rather than performance.
Open implementations of bucket-graph labelling do exist inside full vehicle-routing solvers; for
example BALDES~\cite{baldes}, a branch-cut-and-price solver for vehicle routing, embeds the
labelling in the solver, where it is not readily extractable or reusable as a standalone,
extensible SPPRC library. To our knowledge no open, high-performance, and \emph{extensible} bucket-graph SPPRC
\emph{library} is available to the community. This paper fills that gap. We stress at the outset
that the algorithms we implement are published; our contribution is a faithful, performant, and
extensible \emph{library}, together with a reproducible benchmark, not a new algorithm.

% Statement of contribution
The contributions of this paper are:
\begin{itemize}
  \item \textbf{An open, extensible, header-only library.} We organize the implementation
        around a compile-time \emph{resource concept}: a new SPPRC variant is added by
        implementing a fixed seven-function interface, and resources compose into a single
        label state via template metaprogramming, with no runtime dispatch and a state layout
        fixed at compile time (\autoref{sec:design}). The resource interface itself, across-arc
        concatenation, and the ng-path destination-marking rule are due to the Meta-Solver
        work~\cite{sadykov2026metasolver}; our increment is to realize them as a compile-time
        C++ concept, so that these mechanisms compose for an \emph{arbitrary, user-defined} set
        of resources, not only the built-in ones.
  \item \textbf{A faithful, performant bucket-graph implementation.} We implement the
        bucket-graph algorithm~\cite{sadykov2021bucket} with bidirectional labelling and
        across-arc concatenation, bucket fixing and arc elimination from completion bounds, a
        structure-of-arrays label store with SIMD-accelerated dominance, and a pluggable
        parallel executor (\autoref{sec:background}, \autoref{sec:design}).
  \item \textbf{A reproducible benchmark study.} We compare \texttt{bgspprc} head-to-head against
        PathWyse on shared public instances, with open data and scripts and a parity construction
        that makes the linear-programming bound identical on both sides. At that equal bound
        \texttt{bgspprc} is $1.3\times$--$2.35\times$ faster than PathWyse (and $1.3\times$--$2.3\times$
        faster even when run single-threaded) and runs within
        $1.9\times$--$2.4\times$ of parallel pull labelling~\cite{petersen2025pulllabelling}, a
        different labelling technique for the same problem (\autoref{sec:experiments}).
\end{itemize}

The remainder of the paper is organized as follows. \autoref{sec:background} recalls the SPPRC,
the ng-path relaxation, and the bucket-graph labelling algorithm. \autoref{sec:design} presents
the library design, centred on the resource concept and its composition.
\autoref{sec:experiments} reports the computational study: the head-to-head comparison with
PathWyse, a mode and SIMD ablation, and a comparison with parallel pull labelling for context.
\autoref{sec:conclusions} concludes.

\section{Background}
\label{sec:background}

This section fixes notation and recalls the SPPRC and the bucket-graph labelling algorithm at the
level needed to follow the library design. The algorithm and its analysis are due to Sadykov,
Uchoa and Pessoa~\cite{sadykov2021bucket}, and the resource interface it can be cast behind is due
to~\cite{sadykov2026metasolver}; we summarize rather than reproduce them, and refer the reader to
those papers for the details. Our aim is only to surface the few resource-dependent operations the
library factors out.

Labelling algorithms solve the SPPRC by propagating partial paths as states (\emph{labels}) and
discarding any label dominated by another at the same vertex; we refer to the survey of Irnich and
Desaulniers~\cite{irnich2005survey} for a comprehensive treatment. The dynamic-programming
formulation with resource extension functions and dominance traces to the column-generation pricing
of Desrochers, Desrosiers and Solomon~\cite{desrochers1992vrptw}. Because the elementary variant is
strongly NP-hard, early exact methods enforced elementarity directly through node-resource
labelling~\cite{feillet2004espp}, but this proved too costly at scale. Two ideas made the relaxation
practical: decremental state-space relaxation recovers elementarity
incrementally~\cite{righini2008dssr}, and the \emph{ng-path} relaxation of Baldacci, Mingozzi and
Roberti~\cite{baldacci2011ng} forbids only cycles local to a per-vertex neighbourhood, trading a
small bound loss for a far cheaper subproblem; bounded \emph{bidirectional} labelling, introduced by
Righini and Salani~\cite{righini2006bidir}, halves the effective path length by extending forward
and backward labels to a midpoint and joining them. The bucket-graph algorithm of Sadykov, Uchoa
and Pessoa~\cite{sadykov2021bucket} builds on this bidirectional, ng-relaxed labelling and
accelerates it by partitioning labels into buckets along the main resources, using completion bounds
to fix buckets and eliminate arcs between column-generation iterations; for a broad account of this
algorithmic family we refer to the survey of Costa, Contardo and Desaulniers~\cite{costa2019bpc}.

We work on a directed graph $G=(V,A)$ with source $o$ and sink $d$. Each arc $a=(i,j)$ has cost
$c_a$ and, for each resource $r$, a consumption $d_a^r$; each vertex $v$ carries a window
$[\ell_v^r,u_v^r]$. A path is \emph{resource-feasible} if every resource's accumulated consumption
stays within the window at each visited vertex, and the SPPRC seeks a minimum-cost
resource-feasible $o$--$d$ path. In column generation the $c_a$ are arc reduced costs and we seek
negative-cost paths; the structure is unchanged.

\subsection{The SPPRC and the ng-path relaxation}
\label{sec:sppr-ng}

A resource is advanced along a path by a \emph{resource extension function} (REF). For a disposable
resource (time, capacity) with accumulated value $q$, forward extension along $a=(i,j)$ is
\begin{align}
  q' = \max\{\,q + d_a^r,\ \ell_j^r\,\}, \qquad \text{infeasible if } q' > u_j^r,
  \label{eq:ref-fw}
\end{align}
where the $\max$ models waiting to the lower window; backward extension mirrors this with a $\min$.
One or two resources are designated \emph{main} resources and index the buckets introduced below.

Because requiring elementary paths makes the SPPRC strongly NP-hard, labelling algorithms solve the
\emph{ng-path} relaxation of Baldacci, Mingozzi and Roberti~\cite{baldacci2011ng}, which forbids
only cycles local to a per-vertex neighbourhood of size $n_g$. The relaxation's memory is carried as a bit mask in the label state,
in a compact local encoding; larger $n_g$ tightens the bound at the cost of a larger state. We
report $n_g\in\{8,16,24\}$.

\subsection{Labelling and the bucket graph}
\label{sec:labelling}

The algorithm builds partial paths as \emph{labels}. A label
\begin{align}
  L = \big(v(L),\; c(L),\; q(L),\; \pi(L)\big)
  \label{eq:label}
\end{align}
records its endpoint $v(L)$, accumulated reduced cost $c(L)$, the vector $q(L)\in\mathbb{R}^{m}$ of
the $m\le 2$ \emph{main}-resource values, and a tuple $\pi(L)$ of the remaining per-resource states
(the ng memory and rank-1 bits). A label is extended first along an arc and then on arrival at the
new vertex by the REF~\eqref{eq:ref-fw}, the extensions of the individual resources composing into
that of $L$. The algorithm discards a label $L_2$ when another label $L_1$ at the same vertex
\emph{dominates} it: $L_1$ is no worse on every main resource,
\begin{align}
  q_r(L_1) \preceq_r q_r(L_2)\ \text{ for all } r,\qquad
  {\preceq_r} =
  \begin{cases}
    \le & r \text{ disposable, forward,}\\
    \ge & r \text{ disposable, backward,}\\
    =   & r \text{ non-disposable,}
  \end{cases}
  \label{eq:order}
\end{align}
and cheap enough after a resource-dependent penalty,
\begin{align}
  c(L_1) + \delta(L_1,L_2) \le c(L_2),
  \label{eq:dominance}
\end{align}
where $\delta$ is a resource-dependent penalty charged when $L_1$ is worse than $L_2$ on some
dimension. Each resource also exposes a constant $\delta_{\min}$ ($0$ for every resource shipped here)
that drives a cheap pre-check: if $c(L_1)-c(L_2) > -\delta_{\min}$ the full dominance test is skipped,
which only ever retains a label and so is always safe. For an ordinary disposable resource
$\delta\in\{0,+\infty\}$; the rank-1 cut resource instead makes $\delta$ \emph{finite}, so a cheaper
label can fail to dominate a dearer one; we return to this case in \autoref{sec:r1c}.

The bucket-graph algorithm~\cite{sadykov2021bucket} accelerates the search by partitioning the
labels at each vertex into \emph{buckets} indexed by their main-resource values (a one- or
two-dimensional grid) and processing buckets in a topological order. A \emph{completion bound}
$\underline{c}(B)$ (a lower bound on the cost-to-go from bucket $B$ to the terminal) lets the
algorithm prune: once a pricing threshold $\theta$ is known, a bucket $B$ with best-label cost
$c^\star(B)$ is \emph{fixed} when
\begin{align}
  c^\star(B) + \underline{c}(B) \ge \theta,
  \label{eq:fixing}
\end{align}
since none of its labels can then extend to an improving ($<\theta$) path, and arcs that likewise
cannot lie on any improving path are \emph{eliminated}, shrinking the graph between column-generation
iterations. The search runs \emph{bidirectionally}~\cite{righini2006bidir}, extending forward labels
to a midpoint $\mu$ on the first main resource and backward labels beyond it, and recovers full paths
by joining a forward label $L_f$ and a backward label $L_b$ \emph{across an
arc}~\cite{sadykov2026metasolver} $a=(i,j)$; the join is feasible when their states are compatible
and yields reduced cost
\begin{align}
  c(L_f) + c_a + c(L_b) + \gamma(L_f,L_b),
  \label{eq:concat}
\end{align}
with $\gamma$ a resource-dependent concatenation term. \autoref{alg:label} states the per-call
labelling loop as the library implements it, and \autoref{alg:dom} the resource-generic dominance
test it calls. What matters for the library is that the resource-dependent quantities they
invoke (the REF~\eqref{eq:ref-fw}, the dominance penalty $\delta$ and its lower bound
$\delta_{\min}$, and the concatenation term $\gamma$) are exactly the operations factored behind its
resource interface (\autoref{sec:design}).

\begin{algorithm}[t]
\caption{Bidirectional bucket-graph labelling for one pricing call.}
\label{alg:label}
\begin{algorithmic}[1]
\Require bucket graph $G$, stage $s$, midpoint $\mu$, threshold $\theta$
\Ensure source--sink paths with reduced cost ${}<\theta$
\State compute completion bounds; \textbf{fix} buckets and \textbf{eliminate} arcs via \eqref{eq:fixing}
       \Comment{between-iteration pruning}
\State seed a forward label at $o$ and a backward label at $d$
\ForAll{directions $\mathrm{dir}\in\{\textsc{fw},\textsc{bw}\}$ \textbf{in parallel}}
  \ForAll{buckets $B$ on side $\mathrm{dir}$ in topological order}
           \Comment{fw up to $\mu$, bw beyond}
    \ForAll{labels $L\in B$ and arcs $a$ incident to $v(L)$}
      \State $L' \gets \textsc{Extend}(L,a)$
             \Comment{REF along $a$ then at the new vertex; skip if cost ${}=+\infty$}
      \If{$L'$ is not dominated in its bucket (\autoref{alg:dom})}
        \State insert $L'$, discarding the labels it dominates
      \EndIf
    \EndFor
  \EndFor
\EndFor
\ForAll{forward $L_f$ and backward $L_b$ meeting across an arc $a$}
  \If{their states are compatible} record the path of reduced cost \eqref{eq:concat} \EndIf
\EndFor
\State \Return recorded paths with reduced cost ${}<\theta$
\end{algorithmic}
\end{algorithm}

\begin{algorithm}[t]
\caption{Dominance test: does $L_1$ dominate $L_2$? (same vertex, direction $\mathrm{dir}$, stage $s$)}
\label{alg:dom}
\begin{algorithmic}[1]
\ForAll{main resources $r$}
  \If{$q_r(L_1) \not\preceq_r q_r(L_2)$ by \eqref{eq:order}} \State \Return \textsc{false} \EndIf
\EndFor
\If{$s\in\{\textsc{Heur1},\textsc{Heur2}\}$} \State \Return $c(L_1)\le c(L_2)$
     \Comment{cost-only; ng and rank-1 state ignored} \EndIf
\If{$c(L_1)-c(L_2) > -\delta_{\min}$} \State \Return \textsc{false}
     \Comment{cheap pre-check, skips the full state test} \EndIf
\State \Return $c(L_1) + \delta(L_1,L_2) \le c(L_2)$
       \Comment{$\delta=\sum_r \texttt{domination\_cost}_r$}
\end{algorithmic}
\end{algorithm}

\section{Library design}
\label{sec:design}

The engine of \autoref{sec:background} is generic over the resources a problem carries. Rather
than fix a set of resources or dispatch on a runtime type, \texttt{bgspprc} expresses ``a
resource'' as a C++23 \emph{concept} and composes resources at compile time. This section
presents that interface (\autoref{sec:concept}), how resources compose (\autoref{sec:pack}) and
how a problem is solved (\autoref{sec:using}), the built-in resources (\autoref{sec:builtins})
including the rank-1 cut resource whose dual prices bend the usual dominance rule
(\autoref{sec:r1c}), the multi-stage solve (\autoref{sec:stages}), and the engineering choices
that set the library apart from PathWyse (\autoref{sec:engineering}).

\subsection{The resource concept}
\label{sec:concept}

A resource is any type that supplies a state and the seven operations the labelling engine needs.
The concept makes this requirement explicit and checkable at compile time
(\autoref{lst:concept}). Each extension returns a pair \texttt{(new\_state, extra\_cost)}; an
\texttt{extra\_cost} of $+\infty$ signals infeasibility, so feasibility and cost flow through one
return value.

\begin{lstlisting}[caption={The \texttt{Resource} concept (\texttt{include/bgspprc/resource.h}).},
  label=lst:concept,float=t]
template <typename R>
concept Resource = requires(const R& r, Direction dir, Symmetry sym,
                            typename R::State s, typename R::State s2,
                            int arc_id, int vertex) {
  typename R::State;                                          // per-label state
  { r.symmetric() }            -> std::same_as<bool>;         // symmetric labeling?
  { r.init_state(dir) }        -> std::same_as<typename R::State>;
  { r.extend_along_arc(dir, s, arc_id) }                      // REF along an arc
                               -> std::same_as<std::pair<typename R::State, double>>;
  { r.extend_to_vertex(dir, s, vertex) }                      // REF at a vertex
                               -> std::same_as<std::pair<typename R::State, double>>;
  { r.domination_cost(dir, vertex, s, s2) } -> std::same_as<double>;   // delta in (2)
  { r.concatenation_cost(sym, vertex, s, s2) } -> std::same_as<double>;// fw/bw join
  { r.min_domination_cost() } -> std::same_as<double>;        // delta_min pre-check
};
\end{lstlisting}

The split between \texttt{extend\_along\_arc} and \texttt{extend\_to\_vertex} is deliberate: it
separates the change a resource undergoes while traversing an arc from the change it undergoes on
\emph{arriving} at a vertex. For the ng-path resource this is what keeps across-arc concatenation
correct: the remembered set is remapped along the arc, but the arrival vertex is marked only at
\texttt{extend\_to\_vertex}, so a label joined mid-arc is not charged twice for its endpoint. The
\texttt{domination\_cost} function supplies the penalty $\delta$ of \eqref{eq:dominance},
\texttt{concatenation\_cost} the join term of \autoref{sec:labelling}, and
\texttt{min\_domination\_cost} the constant $\delta_{\min}$ used for the dominance pre-check.

\subsection{Composing resources}
\label{sec:pack}

A problem's resources are bundled into a \texttt{ResourcePack<Rs...>}, a variadic template that
holds the resource objects and presents the same seven operations to the engine. Every operation
fans out over the pack with a fold expression over \texttt{std::index\_sequence}: extensions
accumulate the cost and short-circuit to $+\infty$ if any resource is infeasible, while the
domination and pre-check costs are summed. Because the pack is a compile-time tuple, the label
state is \texttt{std::tuple<typename Rs::State...>}, its size is known at compile time, and no
virtual dispatch occurs on the hot path. The empty pack \texttt{ResourcePack<>} is the
no-extra-resource case (a plain resource-constrained shortest path on the main resources alone).

Adding a new SPPRC variant therefore means writing one type that satisfies
\autoref{lst:concept} and naming it in a pack. \autoref{lst:custom} shows a complete capacity
resource that tracks accumulated arc demand: forward labels accumulate demand and become infeasible
once it exceeds the vehicle capacity $Q$, while backward labels carry the remaining capacity and
fail if it drops below zero, the two directions the bidirectional search of
\autoref{sec:labelling} requires. Because the resource is not symmetric (\texttt{symmetric()} returns
\texttt{false}), forward and backward states have different meanings, and the extension, dominance,
and concatenation functions each branch on the direction.

\begin{lstlisting}[caption={A complete user-defined resource (\texttt{examples/custom\_resource.cpp}, abridged).},
  label=lst:custom,float=t]
struct CapacityResource {
  using State = double;                       // fw: demand used; bw: capacity left
  const double* demand;                       // per-arc demand
  double capacity;                            // Q
  bool symmetric() const { return false; }
  State init_state(Direction dir) const {
    return dir == Direction::Forward ? 0.0 : capacity;
  }
  std::pair<State,double> extend_along_arc(Direction dir, State s, int a) const {
    return {dir == Direction::Forward ? s + demand[a] : s - demand[a], 0.0};
  }
  std::pair<State,double> extend_to_vertex(Direction dir, State s, int) const {
    bool bad = dir == Direction::Forward ? s > capacity + EPS : s < -EPS;
    return {s, bad ? INF : 0.0};              // INF == infeasible
  }
  double domination_cost(Direction dir, int, State s1, State s2) const {
    bool worse = dir == Direction::Forward ? s1 > s2 + EPS : s1 < s2 - EPS;
    return worse ? INF : 0.0;                 // a worse state cannot dominate
  }
  double concatenation_cost(Symmetry, int, State s_fw, State s_bw) const {
    return s_fw + (capacity - s_bw) > capacity + EPS ? INF : 0.0;
  }
  double min_domination_cost() const { return 0.0; }
};
static_assert(bgspprc::Resource<CapacityResource>);          // checked at compile time
using Pack = bgspprc::ResourcePack<CapacityResource>;
\end{lstlisting}

\subsection{Solving a problem}
\label{sec:using}

The graph is passed to the solver as a non-owning \texttt{ProblemView}: arrays of arc
endpoints, costs, per-resource consumptions, and per-vertex windows that the caller continues to
own. A solver is then instantiated on a pack and queried for negative-cost paths.
\autoref{lst:usage} shows the whole flow for a plain resource-constrained shortest path (the empty
pack); swapping \texttt{EmptyPack} for the \texttt{Pack} of \autoref{lst:custom} adds the capacity
constraint with no other change. The pricing threshold $\theta$ selects which paths are returned:
a small negative value keeps only improving columns for column generation, while a large value
enumerates all feasible paths.

\begin{lstlisting}[caption={Solving an SPPRC instance (adapted from \texttt{examples/basic\_spprc.cpp}).},
  label=lst:usage,float=t]
using namespace bgspprc;
ProblemView pv;                                  // non-owning view of caller data
pv.n_vertices = n;  pv.source = 0;  pv.sink = n - 1;
pv.n_arcs = m;      pv.arc_from = from;  pv.arc_to = to;  pv.arc_base_cost = cost;
pv.n_resources = 1; pv.arc_resource = arc_res;             // one resource: time
pv.vertex_lb = v_lb; pv.vertex_ub = v_ub;  pv.n_main_resources = 1;

Solver<EmptyPack> solver(pv, EmptyPack{},
                         {.bucket_steps = {5.0, 1.0}, .theta = -1e-6});  // pricing
solver.set_stage(Stage::Exact);
solver.build();
for (const auto& p : solver.solve())             // paths with reduced cost < theta
  std::printf("cost=%.2f  len=%zu\n", p.reduced_cost, p.vertices.size());
\end{lstlisting}

\subsection{Built-in resources}
\label{sec:builtins}

Five resources ship with the library, covering the variants of~\cite{sadykov2026metasolver}
(\autoref{tbl:resources}). The \emph{standard} resource models time windows and capacity; when it
is a main resource the bucket grid handles its extension directly for speed, and it is otherwise a
pack member. The \emph{ng-path} resource implements the elementarity relaxation of
\autoref{sec:sppr-ng} with a compact local bit representation and the destination-marking split
described above. The \emph{rank-1 cut} resource carries up to 64 limited-memory rank-1 cuts in a
bitset; it is the most instructive of the five because its dual prices give it a dominance rule
unlike any ordinary resource, and we treat it separately in \autoref{sec:r1c}. The \emph{cumulative-cost} and
\emph{pickup-and-delivery} resources cover the cumulative-cost VRP and simultaneous
pickup-and-delivery, the former exercising genuinely non-zero arc cost deltas through
\texttt{extend\_along\_arc}.

\begin{table}[t]
  \centering
  \caption{The five built-in resources and the SPPRC variants they express. The user-defined
    resource of \autoref{lst:custom} adds a sixth.}
  \label{tbl:resources}
  \begin{tabular}{lll}
    \toprule
    Resource & State & Variant / role \\
    \midrule
    Standard           & \texttt{double}            & time windows, capacity \\
    Ng-path            & \texttt{uint32\_t}         & elementarity relaxation \\
    Rank-1 cuts        & \texttt{uint64\_t}         & limited-memory rank-1 cuts \\
    Cumulative cost    & two \texttt{double}s       & cumulative-cost VRP \\
    Pickup-delivery    & two \texttt{double}s       & simultaneous pickup-and-delivery \\
    \bottomrule
  \end{tabular}
\end{table}

\subsection{Rank-1 cuts: a resource that bends dominance}
\label{sec:r1c}

Subset-row cuts, introduced by Jepsen et al.~\cite{jepsen2008src} and generalized to the
limited-memory rank-1 cuts~\cite{pecin2017lmr1c} now central to modern branch-cut-and-price for
routing~\cite{pessoa2020vrpsolver}, are a stringent test of whether an interface is really general,
because their dominance does not fit the usual mould. A rank-1 cut
places a coefficient on routes according to how often they visit a base set $C\subseteq V$; the
separated cut has dual price $\sigma\le 0$, so completing a unit of the cut's coefficient adds a
reduced-cost \emph{penalty} $\beta=-\sigma\ge 0$ to a route: a label that has accumulated more
coefficient is, all else equal, \emph{worse}. What is unusual is that under limited memory a partially
accumulated coefficient may be \emph{forgotten} before it completes (below), so a label that is ahead
on a cut is not infeasible, only at risk of a future $+\beta$. The dominance penalty is therefore
\emph{finite and positive} rather than $\{0,+\infty\}$: a label worse on its cut state can still
dominate if it is cheap enough, and a cheaper label can fail to dominate, neither of which arises for
an ordinary disposable resource.

The library expresses each cut as a single bit of a \texttt{uint64\_t} state (up to 64 active cuts),
specialized to the original three-row, $p=\tfrac12$ subset-row cut. The two extension functions carry the
limited-memory bookkeeping. \texttt{extend\_along\_arc} applies the memory: it ANDs the state with a
per-arc \emph{keep mask}, clearing the bits of any cut whose memory does not contain the arc, so a
route that strays outside a cut's memory forgets its partial accumulation. \texttt{extend\_to\_vertex}
toggles the bits of every cut whose base set contains the vertex; a bit flipping from set to clear
is an overflow (a full coefficient completed) and is charged $+\beta\ge 0$, i.e.\ a penalty
(\autoref{lst:r1c}).

\begin{lstlisting}[caption={The rank-1 cut extension functions (\texttt{include/bgspprc/r1c.h}, abridged).},
  label=lst:r1c,float=t]
std::pair<State,double> extend_along_arc(Direction, State s, int arc) const {
  return {s & arc_keep_mask_[arc], 0.0};        // limited memory: forget outside AM
}
std::pair<State,double> extend_to_vertex(Direction, State s, int v) const {
  uint64_t toggle = vertex_toggle_mask_[v];
  uint64_t overflow = s & toggle;               // bits that flip 1 -> 0
  return {s ^ toggle, overflow ? compute_overflow_cost(overflow) : 0.0};  // penalty >= 0
}
\end{lstlisting}

This contribution flows through the interface without special-casing the engine. Holding a cut bit is
a \emph{disadvantage}: the label is one base-set visit from completing a pair and paying the $+\beta$
penalty, so dominance must account for a penalty the dominator will face but the dominated label will
not. If $L_1$ holds a cut bit that $L_2$ lacks, then on some completion $L_1$ incurs a $+\beta$
penalty that $L_2$ escapes, so $L_1$ must be correspondingly cheaper to dominate: the penalty
$\delta$ of \eqref{eq:dominance} is $\sum_{\ell:\,b_\ell(L_1)\wedge\neg b_\ell(L_2)}\beta_\ell\ge 0$,
one $\beta_\ell$ per such cut, which \texttt{domination\_cost} returns. Unlike an ordinary disposable
resource, whose penalty is $0$ or $+\infty$, the rank-1 resource thus contributes a \emph{finite,
strictly positive} $\delta$: a cheaper label can fail to dominate, and a label worse on its cut state
can still dominate if it is cheap enough. The penalty is never negative (it is zero when $L_1$
carries no pending cut bit that $L_2$ lacks), so \texttt{min\_domination\_cost} returns $0$ and the
pre-check of \autoref{sec:labelling} needs no change. Concatenation is symmetric: when a forward and a
backward label share a cut bit, joining them completes that coefficient, so the join is charged the
same $+\beta$ overflow penalty. Nothing in the bucket-graph engine knows about cuts; the resource
supplies seven functions and the rest follows.

\subsection{Multi-stage solve}
\label{sec:stages}

A full exact labelling is wasteful while many negative-cost columns are easy to find, so the solver
escalates automatically through three stages of increasing strength. \emph{Heuristic~1} keeps only
the single cheapest label per bucket (the most aggressive pruning) and ignores the ng and
rank-1 state; \emph{Heuristic~2} applies full bucket dominance but still ignores ng and rank-1;
\emph{Exact} runs the complete dominance of \eqref{eq:dominance} over every resource. (The
main-resource comparison is always enforced; the stages differ in how much of the remaining state
they consult.) The solver advances when the current stage finds no improving columns (or after a
few iterations), so the cheap stages carry the early column-generation iterations and the exact
engine is reserved for the tail. Because a heuristic stage only prunes more aggressively, it can
\emph{miss} improving columns but never wrongly certify that none exist; the exact stage of
\eqref{eq:dominance} is therefore always run before the solver concludes that no negative-reduced-cost
column remains, so the column-generation loop retains its optimality guarantee. Any column a
heuristic stage does return is a genuine negative-reduced-cost path and may enter the restricted
master directly. A separate \emph{Enumerate} mode, invoked explicitly rather than
as a rung of this ladder, disables dominance and enumerates all paths within the optimality gap up
to a label budget, as required by some rank-1 cut separation routines.

\subsection{Engineering and contrast with PathWyse}
\label{sec:engineering}

The library is header-only and depends only on the C++ standard library; a problem is solved by
instantiating \texttt{Solver<Pack, Executor>} on a non-owning view of the graph data. Three
engineering choices matter for performance. First, labels are stored \emph{structure-of-arrays}:
costs, main resources, and ng/rank-1 bitsets live in parallel per-bucket arrays, so a dominance
scan reads contiguous memory with no gather. Second, dominance is \emph{SIMD-accelerated} over
those arrays (the cost and main-resource pre-filter vectorizes over \texttt{double}, and the ng
subset tests over 32-bit lanes), using the \texttt{<experimental/simd>} parallelism TS where
available, with a scalar fallback (its benefit is instance-dependent; see
\autoref{sec:exp-ablation}). Third, parallelism is delegated to a pluggable \texttt{Executor}
concept: the default executor is sequential, while a thread-pool executor (supplied as the
\texttt{Executor} template argument) runs forward and backward labelling concurrently and chunks the
across-arc join with per-chunk buffers to avoid false sharing.

These choices distinguish \texttt{bgspprc} from PathWyse~\cite{salani2024pathwyse}, the closest
open comparator. PathWyse implements standard labelling with decremental state-space relaxation and
adds resources at runtime by subclassing a resource class; \texttt{bgspprc}
implements the bucket-graph algorithm with bidirectional across-arc concatenation, bucket fixing,
and arc elimination, and composes resources at compile time with no dispatch on the hot path. The
next section quantifies what this buys.

\section{Computational study}
\label{sec:experiments}

We evaluate \texttt{bgspprc} on three public instance families. \emph{spprclib} (45 instances) and
\emph{roberti} (31 instances) are capacity-constrained shortest-path pricing subproblems that arise
in VRP column generation, taken from the \texttt{cptp} repository and used in the
capacitated-profitable-tour study of Jepsen et al.~\cite{jepsen2014cptp}: spprclib comprises the
elementary shortest-path-with-capacity instances of Jepsen, Petersen and
Spoorendonk~\cite{jepsen2008spprc}, and roberti was generated for that study by R.~Roberti using
his pricing code from Baldacci, Mingozzi and Roberti~\cite{baldacci2011ng}. \emph{rcspp} is a standalone
resource-constrained shortest-path benchmark generated as VRPTW column-generation pricing
subproblems from the 56 Solomon instances~\cite{solomon1987} and distributed as the
\texttt{rcspp\_dataset} collection~\cite{spoorendonk2025rcspp}. Each instance is solved at
ng-neighbourhood sizes $n_g\in\{8,16,24\}$ with a 120\,s timeout. Instances are fetched from pinned
public sources and the reference optima are committed; the \texttt{bgspprc} and PathWyse results are
reproducible end to end (rebuild and re-run), and every table and figure regenerates from the
committed result files by a single command.

We report the \emph{shifted geometric mean} (SGM) of runtimes,
$\exp\!\big(\overline{\log(t+s)}\big)-s$ with shift $s=1$\,s, alongside the arithmetic mean and the
number of instances solved within the timeout; timeouts are substituted at 120\,s in both means
rather than dropped, a conservative choice that never rewards a solver for failing. The 1\,s shift
keeps resolution in the sub-second regime where many instances sit. All runs were performed on a
single workstation with an AMD Ryzen~9 3950X processor (16 cores, 32 threads) and 128\,GB of RAM,
running Ubuntu~24.04 (Linux kernel~6.17). The library was built with GCC~14 and \texttt{-march=native}
and run in its default SIMD-enabled \texttt{para\_bidir} mode at the hardware-default thread count
(32 threads); PathWyse was built from a pinned upstream commit with the parity patches of
\autoref{sec:exp-pathwyse} applied automatically. We do not profile memory or per-component
time; the focus is end-to-end pricing runtime. Every table and figure below is regenerated by a
single command from the committed result files.

\subsection{Comparison with PathWyse}
\label{sec:exp-pathwyse}

PathWyse and \texttt{bgspprc} solve the same relaxation but tighten it differently by default, so
a like-for-like \emph{runtime} comparison requires that both compute the same \emph{bound}. We
therefore run both in pure-ng mode: two small patches to PathWyse (both committed with the
benchmark) align its ng-set construction with ours and remove its hard-coded two-cycle
elimination, which is an extra tightening the pure ng-relaxation does not include. With these,
every instance on which both solvers finish agrees on the optimal reduced cost up to cost-scale
rounding, so \autoref{tab:pathwyse} compares runtime at a fixed bound.

\begin{table}[t]
  \centering
  \caption{Runtime comparison of \texttt{bgspprc} against patched PathWyse, both in pure-ng mode, on the spprclib and roberti instance sets. Times in seconds; SGM is the shifted geometric mean ($\exp(\overline{\log(t+1)})-1$); timeouts (120\,s) are substituted, never dropped. The bgspprc columns are the default 32-thread \texttt{para\_bidir} build (sgm, mean, solved) and the single-threaded \texttt{bidir\_base} build ($\text{sgm}_1$); both speedup columns are $(\text{pw}_{\text{sgm}}+1)/(\text{bg}_{\text{sgm}}+1)$ against PathWyse, for the 32-thread and 1-thread builds respectively. bgspprc is faster than single-threaded PathWyse even when itself run single-threaded.}
  \label{tab:pathwyse}
  \small
  \begin{tabular}{llrrrrrrrrr}
    \toprule
    Set & $n_g$ & \multicolumn{4}{c}{bgspprc} & \multicolumn{3}{c}{PathWyse} & \multicolumn{2}{c}{Speedup} \\
    \midrule
     &  & sgm & mean & solved & $\text{sgm}_1$ & sgm & mean & solved & 32-thr & 1-thr \\
    \midrule
    spprclib & 8 & 0.917 & 7.238 & 44/45 & 0.922 & 1.522 & 12.133 & 42/45 & 1.32$\times$ & 1.31$\times$ \\
    spprclib & 16 & 2.048 & 15.801 & 40/45 & 1.997 & 4.284 & 24.290 & 38/45 & 1.73$\times$ & 1.76$\times$ \\
    spprclib & 24 & 5.838 & 26.725 & 38/45 & 5.286 & 10.590 & 42.512 & 31/45 & 1.69$\times$ & 1.84$\times$ \\
    roberti & 8 & 0.551 & 0.911 & 31/31 & 0.636 & 2.330 & 9.999 & 30/31 & 2.15$\times$ & 2.03$\times$ \\
    roberti & 16 & 3.176 & 15.003 & 28/31 & 3.239 & 8.796 & 28.152 & 28/31 & 2.35$\times$ & 2.31$\times$ \\
    roberti & 24 & 14.813 & 44.346 & 23/31 & 14.247 & 26.824 & 62.019 & 19/31 & 1.76$\times$ & 1.82$\times$ \\
    \bottomrule
  \end{tabular}
\end{table}

\begin{figure}[t]
  \centering
  \includegraphics[width=0.7\linewidth]{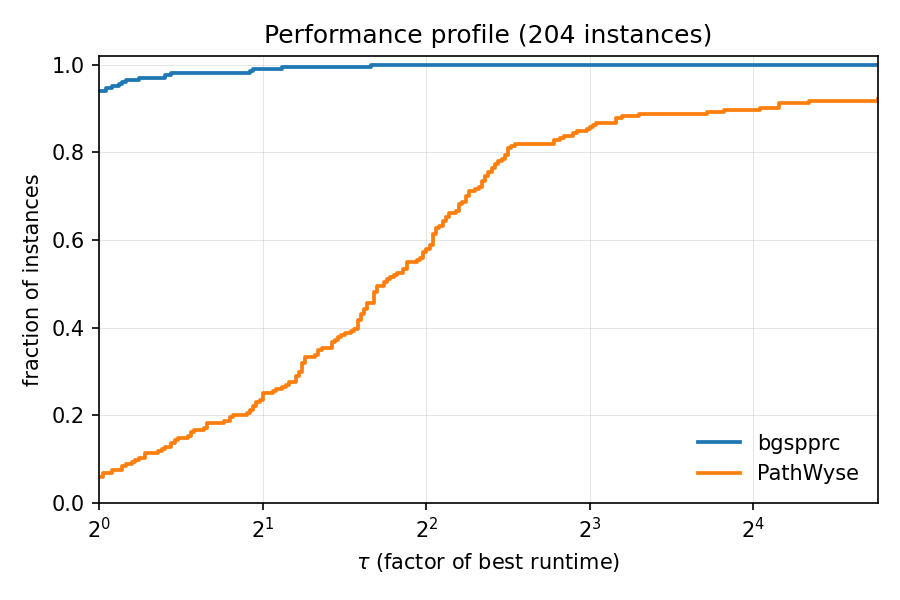}
  \caption{Performance profile on the shared spprclib and roberti instances: for each solver,
    the fraction of instances (vertical axis) solved within a factor $\tau$ of the fastest
    solver's time on that instance (horizontal axis, log scale). A higher curve is better.}
  \label{fig:perf-profile}
\end{figure}

Across all six (set, $n_g$) cells, \texttt{bgspprc} is faster than PathWyse by shifted geometric
mean, with per-cell speedups of $1.3\times$--$2.35\times$ (\autoref{tab:pathwyse}); it also solves
at least as many instances within the timeout in every cell, and strictly more in five of the
six (tying only on roberti at $n_g=16$). The advantage is largest on the roberti instances
at $n_g=8$ and $n_g=16$ ($2.15\times$ and $2.35\times$). The performance profile (\autoref{fig:perf-profile}) shows the same conclusion
instance by instance: \texttt{bgspprc} is the faster solver on the large majority of instances and
its curve dominates PathWyse's throughout. Since both solvers compute the same bound, this is a
runtime advantage at equal strength, which we attribute to the bucket-graph pruning and the
structure-of-arrays label store of \autoref{sec:engineering}.

Two further observations support reading these numbers as a genuine result at equal bound. First,
on correctness: of the 228 (instance, $n_g$) pairs, both solvers finish within the timeout on 188.
On 152 of those the optimal reduced cost is bit-identical and on 185 the two agree to within
$0.002$, a rounding artifact of PathWyse's integer cost scaling. Of the three pairs that differ by
more, \texttt{bgspprc} reports a strictly \emph{better} bound on two (by up to $0.137$) and a bound
worse by only $0.002$ (itself at the rounding tolerance) on the third; \texttt{bgspprc}'s optimum
is thus never more than $0.002$ worse than PathWyse's, confirming that the parity construction yields
the same LP bound and that our implementation is correct against an independent reference. Second, on
the head-to-head: among the 188 jointly solved pairs \texttt{bgspprc} is strictly faster on 176 and
slower on only 12, so the SGM advantage is not an artifact of a few large instances but holds
broadly across the benchmark.

A note on compute. \autoref{tab:pathwyse} reports each solver in its default configuration:
\texttt{bgspprc}'s \texttt{para\_bidir} mode runs forward and backward labelling concurrently
(32 threads here), whereas PathWyse is single-threaded by design. We compare the two libraries as
each is intended to be run; exploiting parallelism is part of the implementation, not a thumb on the
scale. That the advantage is nonetheless \emph{algorithmic} rather than a parallelism artifact is
shown by the $\text{sgm}_1$ and 1-thr columns of \autoref{tab:pathwyse}, which repeat the comparison
with \texttt{bgspprc} restricted to a single thread (its sequential \texttt{bidir\_base} build): even
then it is $1.3\times$--$2.3\times$ faster than PathWyse at the identical bound in every cell and
solves at least as many instances. At the largest neighbourhoods ($n_g=24$) the single-threaded build
is in fact faster than the parallel one, whose thread-pool overhead does not pay off there. We
therefore attribute the advantage primarily to the bucket-graph pruning and the structure-of-arrays
label store of \autoref{sec:engineering}, with thread-level parallelism adding a further margin on
top.

The gap between the SGM and the arithmetic mean in \autoref{tab:pathwyse} (the arithmetic mean
is several times larger for both solvers) reflects a heavy tail of hard instances near the
timeout that both solvers share.

\subsection{Mode and SIMD ablation}
\label{sec:exp-ablation}

To see where the runtime goes, \autoref{tab:ablation} ablates the two execution axes on the rcspp
set: monodirectional versus bidirectional versus parallel-bidirectional search, each with SIMD on
or off. Two patterns stand out. Monodirectional search is fastest at $n_g=8$, where the search
trees are shallow and the bookkeeping of splitting into forward and backward halves does not pay
off; from $n_g=16$ upward the parallel-bidirectional mode catches up and stays within a sub-second
SGM of the leader. The one outlier is \texttt{para\_bidir\_base} at $n_g=8$ (SGM 8.80\,s, 37/56
solved), a real cliff for the scalar parallel-bidirectional mode on the shallowest setting that the
SIMD default (\texttt{para\_bidir\_vec}, 1.88\,s, 56/56) avoids.

The SIMD axis deserves comment, because vectorizing the dominance check does \emph{not} help here:
for the mono and bidir modes \texttt{\_vec} and \texttt{\_base} differ by under 5\% on rcspp (the
\texttt{para\_bidir\_base} cliff at $n_g=8$ aside), and on the larger spprclib and roberti instances
(\autoref{tab:ablation-full}) the vectorized build is frequently \emph{slower}, by up to about
$40\%$. This follows from how the two scan a bucket. The scalar check exploits the
cost-sorted label store: it begins at the first label that could be dominated and stops at the
first dominator, so it usually settles a label after a handful of comparisons. By contrast, the
vectorized check evaluates a fixed-width batch of candidates against each stored label before it can
stop. The batch only amortizes its overhead when many labels per bucket survive the cost and
main-resource pre-filter and must be compared in full; on these pure-ng instances, with one or two
main resources and cost-sorted buckets, that surviving set is small and the scalar early-exit wins.
We therefore \emph{conjecture} that SIMD dominance pays off precisely when this early termination is
defeated: when the non-dominated frontier per bucket is wide and the cost pre-filter is weak. The
clearest such regime is pricing with active limited-memory rank-1 cuts (\autoref{sec:r1c}): there a
cut's dominance penalty $\delta$ is finite and positive rather than $\{0,+\infty\}$, so a cheaper
label can fail to dominate while a label worse on its cut state can still dominate. Both effects
defeat the cost-sorted early-exit and enlarge the set of mutually non-dominated labels, exactly the
conditions the vectorized scan is built for. Our benchmark runs pure ng with no cuts, the regime in
which SIMD has least to do; evaluating it under rank-1-cut and higher-dimensional resource pricing
is left as future work. One consequence is that the headline comparison is \emph{conservative}: it
reports the SIMD default (\texttt{para\_bidir\_vec}) throughout, even though the scalar configurations are
often faster still.

\autoref{tab:axis-summary} distils the six modes into one SGM speedup ratio per execution axis,
each isolated on the scalar, single-thread, bidirectional baseline, and confirms three patterns.
\emph{Direction} is the dominant axis and changes sign with search depth: monodirectional search is
faster on the shallow rcspp instances and on spprclib at $n_g=8$ ($0.4$--$0.6\times$), while
bidirectional search is up to $3.1\times$ faster once the search deepens (spprclib $n_g\ge16$,
roberti throughout). \emph{Threads} add a smaller, near-constant $1.02$--$1.27\times$ on top of
bidirectional search everywhere except the scalar rcspp $n_g=8$ cliff. \emph{SIMD} never helps on
these pure-ng instances ($0.74$--$0.99\times$); its only role is to rescue that cliff. The fastest
build in all but one cell is therefore a scalar one (monodirectional when the search is
shallow, scalar parallel-bidirectional when it is deep; the lone exception is rcspp $n_g=8$, where
vectorized monodirectional labelling wins by a negligible margin), yet each of those has a regime in which it
collapses: \texttt{para\_bidir\_base} times out on 19 of the 56 easiest rcspp $n_g=8$ instances,
and \texttt{mono} is roughly $3\times$ slower and solves fewer of the deep roberti instances. The
shipped default \texttt{para\_bidir\_vec} is never the per-cell fastest, but it is the only build
with no catastrophic regime, robust rather than universally optimal. A search-depth-adaptive policy
that selects monodirectional labelling for the smallest subproblems and parallel-bidirectional
labelling otherwise would capture the per-cell best of both; we leave it to future work.

\begin{table}[t]
  \centering
  \caption{Mode and SIMD ablation of \texttt{bgspprc} on the 56 Solomon RCSPP instances, by $n_g$. The bidirectional axis (\texttt{mono} / \texttt{bidir} / \texttt{para\_bidir}) crosses a SIMD axis (\texttt{\_base} = scalar, \texttt{\_vec} = SIMD). Times in seconds; same SGM and 120\,s timeout convention as Table~\ref{tab:pathwyse}. The default CLI build is \texttt{para\_bidir\_vec} (invoked as \texttt{para\_bidir}). The full table over all three instance sets is Table~\ref{tab:ablation-full}.}
  \label{tab:ablation}
  \begin{tabular}{llrrr}
    \toprule
    $n_g$ & mode & sgm (s) & mean (s) & solved \\
    \midrule
    8 & \texttt{mono\_base} & 1.296 & 2.885 & 56/56 \\
     & \texttt{mono\_vec} & 1.288 & 2.850 & 56/56 \\
     & \texttt{bidir\_base} & 2.389 & 7.870 & 56/56 \\
     & \texttt{bidir\_vec} & 2.442 & 7.995 & 56/56 \\
     & \texttt{para\_bidir\_base} & 8.798 & 44.624 & 37/56 \\
     & \texttt{para\_bidir\_vec} & 1.884 & 5.640 & 56/56 \\
    \midrule
    16 & \texttt{mono\_base} & 1.578 & 4.979 & 56/56 \\
     & \texttt{mono\_vec} & 1.598 & 5.134 & 56/56 \\
     & \texttt{bidir\_base} & 2.824 & 12.800 & 55/56 \\
     & \texttt{bidir\_vec} & 2.855 & 12.952 & 55/56 \\
     & \texttt{para\_bidir\_base} & 2.222 & 8.913 & 56/56 \\
     & \texttt{para\_bidir\_vec} & 2.200 & 8.783 & 56/56 \\
    \midrule
    24 & \texttt{mono\_base} & 1.940 & 11.005 & 54/56 \\
     & \texttt{mono\_vec} & 1.974 & 11.567 & 54/56 \\
     & \texttt{bidir\_base} & 3.156 & 17.821 & 51/56 \\
     & \texttt{bidir\_vec} & 3.222 & 18.443 & 51/56 \\
     & \texttt{para\_bidir\_base} & 2.532 & 14.792 & 52/56 \\
     & \texttt{para\_bidir\_vec} & 2.588 & 15.433 & 52/56 \\
    \bottomrule
  \end{tabular}
\end{table}

\begin{table}[t]
  \centering
  \caption{Per-axis SGM speedup decomposition of the mode/SIMD ablation (Tables~\ref{tab:ablation} and~\ref{tab:ablation-full}). Each ratio is $\text{sgm}_{\text{lighter}}/\text{sgm}_{\text{heavier}}$, so a value $>1$ means the heavier option (bidirectional, parallel, or vectorized) is faster; each axis is isolated on the scalar, single-thread, bidirectional baseline. \emph{Direction} $=$ \texttt{mono\_base}/\texttt{bidir\_base}; \emph{Threads} $=$ \texttt{bidir\_base}/\texttt{para\_bidir\_base}; \emph{SIMD} $=$ \texttt{bidir\_base}/\texttt{bidir\_vec}. The \emph{fastest} column is the minimum-SGM build of all six. Direction changes sign with search depth, threading adds a near-constant boost (bar the scalar rcspp $n_g=8$ cliff), and SIMD never helps on these pure-ng instances; the default \texttt{para\_bidir\_vec} is never the per-cell fastest but is the only build with no catastrophic regime.}
  \label{tab:axis-summary}
  \begin{tabular}{llrrrl}
    \toprule
    set & $n_g$ & Direction & Threads & SIMD & fastest \\
    \midrule
     &  & mono$\to$bidir & bidir$\to$para & base$\to$vec & mode \\
    \midrule
    rcspp & 8 & 0.54$\times$ & 0.27$\times$ & 0.98$\times$ & \texttt{mono\_vec} \\
     & 16 & 0.56$\times$ & 1.27$\times$ & 0.99$\times$ & \texttt{mono\_base} \\
     & 24 & 0.61$\times$ & 1.25$\times$ & 0.98$\times$ & \texttt{mono\_base} \\
    \midrule
    spprclib & 8 & 0.42$\times$ & 1.02$\times$ & 0.98$\times$ & \texttt{mono\_base} \\
     & 16 & 1.18$\times$ & 1.06$\times$ & 0.90$\times$ & \texttt{para\_bidir\_base} \\
     & 24 & 3.09$\times$ & 1.09$\times$ & 0.81$\times$ & \texttt{para\_bidir\_base} \\
    \midrule
    roberti & 8 & 1.38$\times$ & 1.25$\times$ & 0.91$\times$ & \texttt{para\_bidir\_base} \\
     & 16 & 2.39$\times$ & 1.22$\times$ & 0.81$\times$ & \texttt{para\_bidir\_base} \\
     & 24 & 3.08$\times$ & 1.26$\times$ & 0.74$\times$ & \texttt{para\_bidir\_base} \\
    \bottomrule
  \end{tabular}
\end{table}

\subsection{Context: comparison with parallel pull labelling}
\label{sec:exp-rcspp}

Finally, we place \texttt{bgspprc} against the parallel pull-labelling algorithm of Petersen and
Spoorendonk~\cite{petersen2025pulllabelling}, as implemented in the Flowty solver~\cite{flowty}, on
the rcspp set (\autoref{tab:rcspp}). The pull-labelling runtimes are those reported in that paper on
the same instances (committed as \texttt{pull\_algo\_runtimes.csv}). Pull labelling
is a different, independently developed labelling technique for the same problem (in principle
equally general) rather than a specialized single-purpose solver; the two are distinct points in
the algorithmic design space. \texttt{bgspprc} runs within $1.9\times$--$2.4\times$ of it in SGM on
the same instances ($n_g\in\{8,16\}$; one fewer at $n_g=24$, 52 versus 53 of 56), and the gap
narrows monotonically as $n_g$ grows ($2.40\times \to 2.10\times \to 1.92\times$), suggesting a
fixed per-instance overhead rather than one that scales with the search. The comparison is best read
as locating \texttt{bgspprc} against a strong, independent alternative on shared instances rather than
as a generality-versus-specialization trade-off; that \texttt{bgspprc} stays within a small constant
factor while adding compile-time extensibility and built-in cut support is the takeaway. We do not
measure how the two methods compare across graph densities (a pull step that gathers each label from
its predecessors has less to amortize on the sparse pricing graphs used here than on dense ones), but
this is conjecture, and a characterization across densities is left to future work.
\begin{table}[t]
  \centering
  \caption{Runtime of \texttt{bgspprc} (\texttt{para\_bidir}) against the parallel pull-labelling algorithm of Petersen \& Spoorendonk (2025) on 56 Solomon RCSPP instances per $n_g$. Times in seconds, same SGM and 120\,s timeout convention as Table~\ref{tab:pathwyse}. Slowdown $=(\text{bg}_{\text{sgm}}+1)/(\text{pull}_{\text{sgm}}+1)$; pull labelling, a different technique, is faster on these instances.}
  \label{tab:rcspp}
  \begin{tabular}{lrrrrrrr}
    \toprule
    $n_g$ & \multicolumn{3}{c}{bgspprc} & \multicolumn{3}{c}{pull-labelling} & Slowdown \\
    \midrule
     & sgm & mean & solved & sgm & mean & solved &  \\
    \midrule
    8 & 1.884 & 5.640 & 56/56 & 0.203 & 0.406 & 56/56 & 2.40$\times$ \\
    16 & 2.200 & 8.783 & 56/56 & 0.526 & 3.010 & 56/56 & 2.10$\times$ \\
    24 & 2.588 & 15.433 & 52/56 & 0.873 & 9.123 & 53/56 & 1.92$\times$ \\
    \bottomrule
  \end{tabular}
\end{table}

\section{Conclusions}
\label{sec:conclusions}

We have presented \texttt{bucket-graph-spprc}, an open-source, header-only C++23 library for the
SPPRC. It implements the bucket-graph labelling algorithm with bidirectional across-arc
concatenation, bucket fixing, and arc elimination, and is organized around a compile-time resource
concept that lets a new SPPRC variant be added by implementing a fixed seven-function interface,
with resources composed into the label state with no runtime dispatch. On shared public instances, at
an equal LP bound, \texttt{bgspprc} is $1.3\times$--$2.35\times$ faster than PathWyse, the main
open-source comparator, and runs within $1.9\times$--$2.4\times$ of parallel pull labelling, a
different labelling technique for the same problem, while remaining general and extensible. The library, the benchmark scripts, and the pinned instances are publicly available
under a permissive licence, and every table and figure in this paper is reproduced by a single
command.

The algorithms are not new; the contribution is an open, performant, and extensible
implementation, together with a reproducible benchmark, of methods that were previously available
only in academic or closed solvers. Several directions remain. The resource interface invites
further built-in resources beyond the five shipped here; integrating the pricer into a full
branch-cut-and-price loop would exercise it under repeated re-solves with cuts and branching; and
the gap to parallel pull labelling, which our ablation suggests is a fixed per-instance
overhead, is worth closing. The library is intended to lower the barrier to building on the bucket-graph
algorithm.

\section*{Code and data availability}

The \texttt{bucket-graph-spprc} library is open source under the MIT licence at
\url{https://github.com/spoorendonk/bucket-graph-spprc}. Its \texttt{benchmarks/} directory
contains the run and comparison scripts, the two PathWyse parity patches, the reference optima, and
the committed result files from which every table and figure in this paper is regenerated by a
single command. The instance sets are fetched from public repositories: the spprclib and roberti
pricing instances from \url{https://github.com/spoorendonk/cptp} and the Solomon RCSPP set from
\url{https://github.com/spoorendonk/rcspp_dataset}. PathWyse is cloned from its upstream repository
and built from a pinned commit with the two parity patches applied automatically. An archival
snapshot of the library (v0.1.0) is deposited on Zenodo~\cite{bgspprc_software}: the version DOI
\texttt{10.5281/zenodo.20819209} pins the exact release this paper describes, while the concept DOI
\texttt{10.5281/zenodo.20819208} always resolves to the latest version.

\section*{Disclosure statement}

The author is affiliated with Flowty (\url{https://flowty.ai}), whose solver provides the parallel
pull-labelling implementation~\cite{flowty} used as a comparator in \autoref{sec:exp-rcspp}. All
benchmark scripts, instances, and reference results are public, so the comparison is independently
reproducible.

\appendix

\section{Full mode and SIMD ablation}
\label{sec:appendix-ablation}

\autoref{tab:ablation-full} gives the mode and SIMD ablation of \autoref{sec:exp-ablation} across
all three instance sets, complementing the rcspp slice in \autoref{tab:ablation}. The pattern is
consistent with the main text. On the smaller spprclib instances monodirectional search is fastest at
$n_g=8$, but the bidirectional modes overtake it from $n_g=16$ onward as the search deepens. On
roberti the parallel-bidirectional mode is fastest at every $n_g$, in its scalar form
\texttt{para\_bidir\_base} throughout. Across both sets the \texttt{\_vec} builds
run slower than their \texttt{\_base} counterparts, by up to about
$40\%$ on the largest instances, so the fastest configuration is frequently a scalar one;
\autoref{sec:exp-ablation} explains the mechanism and conjectures when vectorization would instead
pay off.

\begin{table}[t]
  \centering
  \caption{Mode and SIMD ablation of \texttt{bgspprc} across all three instance sets, by $n_g$ (companion to Table~\ref{tab:ablation}). Same conventions.}
  \label{tab:ablation-full}
  \footnotesize
  \begin{tabular}{lllrrr}
    \toprule
    set & $n_g$ & mode & sgm (s) & mean (s) & solved \\
    \midrule
    spprclib & 8 & \texttt{mono\_base} & 0.390 & 3.026 & 44/45 \\
     &  & \texttt{mono\_vec} & 0.429 & 3.094 & 44/45 \\
     &  & \texttt{bidir\_base} & 0.922 & 7.170 & 44/45 \\
     &  & \texttt{bidir\_vec} & 0.941 & 7.309 & 44/45 \\
     &  & \texttt{para\_bidir\_base} & 0.905 & 7.169 & 44/45 \\
     &  & \texttt{para\_bidir\_vec} & 0.917 & 7.238 & 44/45 \\
    \midrule
    spprclib & 16 & \texttt{mono\_base} & 2.354 & 7.058 & 44/45 \\
     &  & \texttt{mono\_vec} & 3.336 & 9.628 & 44/45 \\
     &  & \texttt{bidir\_base} & 1.997 & 15.773 & 40/45 \\
     &  & \texttt{bidir\_vec} & 2.224 & 16.100 & 40/45 \\
     &  & \texttt{para\_bidir\_base} & 1.875 & 15.586 & 40/45 \\
     &  & \texttt{para\_bidir\_vec} & 2.048 & 15.801 & 40/45 \\
    \midrule
    spprclib & 24 & \texttt{mono\_base} & 16.323 & 46.133 & 34/45 \\
     &  & \texttt{mono\_vec} & 22.270 & 55.078 & 30/45 \\
     &  & \texttt{bidir\_base} & 5.286 & 25.308 & 38/45 \\
     &  & \texttt{bidir\_vec} & 6.533 & 27.552 & 38/45 \\
     &  & \texttt{para\_bidir\_base} & 4.869 & 24.788 & 38/45 \\
     &  & \texttt{para\_bidir\_vec} & 5.838 & 26.725 & 38/45 \\
    \midrule
    roberti & 8 & \texttt{mono\_base} & 0.880 & 1.470 & 31/31 \\
     &  & \texttt{mono\_vec} & 1.022 & 1.774 & 31/31 \\
     &  & \texttt{bidir\_base} & 0.636 & 1.133 & 31/31 \\
     &  & \texttt{bidir\_vec} & 0.697 & 1.246 & 31/31 \\
     &  & \texttt{para\_bidir\_base} & 0.511 & 0.851 & 31/31 \\
     &  & \texttt{para\_bidir\_vec} & 0.551 & 0.911 & 31/31 \\
    \midrule
    roberti & 16 & \texttt{mono\_base} & 7.750 & 28.426 & 27/31 \\
     &  & \texttt{mono\_vec} & 10.484 & 32.769 & 25/31 \\
     &  & \texttt{bidir\_base} & 3.239 & 15.128 & 28/31 \\
     &  & \texttt{bidir\_vec} & 3.981 & 16.270 & 28/31 \\
     &  & \texttt{para\_bidir\_base} & 2.649 & 14.165 & 28/31 \\
     &  & \texttt{para\_bidir\_vec} & 3.176 & 15.003 & 28/31 \\
    \midrule
    roberti & 24 & \texttt{mono\_base} & 43.841 & 70.594 & 17/31 \\
     &  & \texttt{mono\_vec} & 57.345 & 79.861 & 15/31 \\
     &  & \texttt{bidir\_base} & 14.247 & 42.547 & 24/31 \\
     &  & \texttt{bidir\_vec} & 19.242 & 48.893 & 23/31 \\
     &  & \texttt{para\_bidir\_base} & 11.300 & 38.679 & 25/31 \\
     &  & \texttt{para\_bidir\_vec} & 14.813 & 44.346 & 23/31 \\
    \midrule
    rcspp & 8 & \texttt{mono\_base} & 1.296 & 2.885 & 56/56 \\
     &  & \texttt{mono\_vec} & 1.288 & 2.850 & 56/56 \\
     &  & \texttt{bidir\_base} & 2.389 & 7.870 & 56/56 \\
     &  & \texttt{bidir\_vec} & 2.442 & 7.995 & 56/56 \\
     &  & \texttt{para\_bidir\_base} & 8.798 & 44.624 & 37/56 \\
     &  & \texttt{para\_bidir\_vec} & 1.884 & 5.640 & 56/56 \\
    \midrule
    rcspp & 16 & \texttt{mono\_base} & 1.578 & 4.979 & 56/56 \\
     &  & \texttt{mono\_vec} & 1.598 & 5.134 & 56/56 \\
     &  & \texttt{bidir\_base} & 2.824 & 12.800 & 55/56 \\
     &  & \texttt{bidir\_vec} & 2.855 & 12.952 & 55/56 \\
     &  & \texttt{para\_bidir\_base} & 2.222 & 8.913 & 56/56 \\
     &  & \texttt{para\_bidir\_vec} & 2.200 & 8.783 & 56/56 \\
    \midrule
    rcspp & 24 & \texttt{mono\_base} & 1.940 & 11.005 & 54/56 \\
     &  & \texttt{mono\_vec} & 1.974 & 11.567 & 54/56 \\
     &  & \texttt{bidir\_base} & 3.156 & 17.821 & 51/56 \\
     &  & \texttt{bidir\_vec} & 3.222 & 18.443 & 51/56 \\
     &  & \texttt{para\_bidir\_base} & 2.532 & 14.792 & 52/56 \\
     &  & \texttt{para\_bidir\_vec} & 2.588 & 15.433 & 52/56 \\
    \bottomrule
  \end{tabular}
\end{table}

\section{Per-instance results}
\label{sec:appendix-per-instance}

For completeness and to back the aggregated cells of \autoref{tab:pathwyse} and
\autoref{tab:rcspp}, this appendix lists every instance's wall-clock time. \autoref{tab:per-instance-pw}
covers the spprclib and roberti sets (\texttt{bgspprc} \texttt{para\_bidir} versus patched PathWyse,
both pure-ng); \autoref{tab:per-instance-pull} covers the rcspp set (\texttt{bgspprc} versus parallel
pull labelling). Times are in seconds at each $n_g\in\{8,16,24\}$, with \textit{t/o} marking a 120\,s
timeout; the shifted geometric mean of each solver column reproduces the corresponding cell of the
in-text tables. These tables are regenerated from the committed result files by the same single
command as the rest of the paper.

\begin{footnotesize}
\begin{longtable}{lrrrrrr}
  \caption{Per-instance wall-clock time (s) on the spprclib and roberti sets: \texttt{bgspprc} (\texttt{para\_bidir}, column \emph{bg}) against patched PathWyse (\emph{pw}), both in pure-ng mode. \textit{t/o} marks a 120\,s timeout. Per-column shifted geometric means within each set block reproduce \autoref{tab:pathwyse}.}\label{tab:per-instance-pw}\\
  \toprule
  Instance & \multicolumn{2}{c}{$n_g=8$} & \multicolumn{2}{c}{$n_g=16$} & \multicolumn{2}{c}{$n_g=24$} \\
   & bg & pw & bg & pw & bg & pw \\
  \midrule
  \endfirsthead
  \multicolumn{7}{l}{\small\itshape (\autoref{tab:per-instance-pw} continued)}\\
  \toprule
  Instance & \multicolumn{2}{c}{$n_g=8$} & \multicolumn{2}{c}{$n_g=16$} & \multicolumn{2}{c}{$n_g=24$} \\
   & bg & pw & bg & pw & bg & pw \\
  \midrule
  \endhead
  \midrule \multicolumn{7}{r}{\small\itshape continued on next page}\\
  \endfoot
  \bottomrule
  \endlastfoot
  \multicolumn{7}{l}{\textbf{spprclib}} \\
  A-n54-k7-149 & 0.051 & 0.153 & 0.311 & 0.674 & 2.372 & 2.726 \\
  A-n60-k9-57 & 0.063 & 0.147 & 0.338 & 0.528 & 1.705 & 1.464 \\
  A-n61-k9-80 & 0.067 & 0.155 & 0.270 & 0.869 & 2.220 & 4.204 \\
  A-n62-k8-99 & 0.114 & 0.533 & 0.915 & 3.342 & 8.839 & 32.199 \\
  A-n63-k10-44 & 0.049 & 0.131 & 0.237 & 0.538 & 1.192 & 1.546 \\
  A-n63-k9-157 & 0.048 & 0.152 & 0.196 & 0.582 & 0.884 & 1.751 \\
  A-n64-k9-45 & 0.098 & 0.423 & 0.349 & 1.810 & 1.630 & 8.344 \\
  A-n65-k9-10 & 0.083 & 0.196 & 0.367 & 0.975 & 2.030 & 3.544 \\
  A-n69-k9-42 & 0.073 & 0.227 & 0.387 & 0.897 & 1.716 & 2.506 \\
  A-n80-k10-14 & 0.771 & 2.043 & 5.698 & 17.017 & 23.151 & \textit{t/o} \\
  B-n45-k6-54 & 0.062 & 0.253 & 0.082 & 0.558 & 0.325 & 6.576 \\
  B-n50-k8-40 & 0.044 & 0.169 & 0.427 & 1.329 & 4.526 & 7.872 \\
  B-n52-k7-15 & 0.154 & 0.226 & 5.355 & 22.080 & 39.616 & \textit{t/o} \\
  B-n57-k7-20 & 0.537 & 4.895 & 0.613 & 5.465 & 0.652 & 10.659 \\
  B-n66-k9-50 & 0.136 & 1.329 & 5.625 & 39.999 & 38.130 & \textit{t/o} \\
  B-n67-k10-26 & 0.093 & 0.371 & 0.948 & 1.836 & 9.566 & 9.693 \\
  B-n68-k9-65 & 0.158 & 0.501 & 1.112 & 6.245 & 31.232 & \textit{t/o} \\
  B-n78-k10-70 & 0.171 & 0.507 & 1.432 & 5.158 & 103.154 & \textit{t/o} \\
  E-n101-k14-158 & 0.129 & 0.418 & 0.695 & 1.377 & 4.428 & 4.435 \\
  E-n101-k8-291 & 0.194 & 0.993 & 1.312 & 5.261 & 11.985 & 22.989 \\
  E-n76-k10-72 & 0.129 & 0.503 & 0.880 & 3.786 & 7.413 & 33.895 \\
  E-n76-k14-102 & 0.064 & 0.089 & 0.179 & 0.241 & 0.679 & 0.648 \\
  E-n76-k15-40 & 0.036 & 0.080 & 0.177 & 0.212 & 0.538 & 0.492 \\
  E-n76-k7-44 & 0.214 & 0.660 & 1.633 & 5.286 & 17.856 & 60.875 \\
  G-n262-k25-316 & \textit{t/o} & \textit{t/o} & \textit{t/o} & \textit{t/o} & \textit{t/o} & \textit{t/o} \\
  M-n101-k10-97 & 0.365 & 2.097 & 2.194 & 28.595 & 33.746 & \textit{t/o} \\
  M-n121-k7-260 & 43.552 & \textit{t/o} & \textit{t/o} & \textit{t/o} & \textit{t/o} & \textit{t/o} \\
  M-n151-k12-15 & 27.692 & 31.453 & \textit{t/o} & \textit{t/o} & \textit{t/o} & \textit{t/o} \\
  M-n200-k16-143 & 5.308 & 94.768 & 6.481 & \textit{t/o} & 7.926 & \textit{t/o} \\
  M-n200-k17-12 & 65.031 & \textit{t/o} & \textit{t/o} & \textit{t/o} & \textit{t/o} & \textit{t/o} \\
  P-n101-k4-174 & 1.988 & 8.695 & 15.499 & \textit{t/o} & \textit{t/o} & \textit{t/o} \\
  P-n50-k10-24 & 0.014 & 0.019 & 0.027 & 0.039 & 0.041 & 0.056 \\
  P-n50-k7-92 & 0.018 & 0.053 & 0.039 & 0.179 & 0.065 & 0.361 \\
  P-n50-k8-19 & 0.052 & 0.114 & 0.096 & 0.178 & 0.183 & 0.512 \\
  P-n51-k10-30 & 0.018 & 0.023 & 0.060 & 0.054 & 0.170 & 0.090 \\
  P-n55-k10-44 & 0.016 & 0.033 & 0.037 & 0.089 & 0.063 & 0.169 \\
  P-n55-k15-88 & 0.018 & 0.020 & 0.030 & 0.023 & 0.060 & 0.028 \\
  P-n55-k7-116 & 0.037 & 0.123 & 0.144 & 0.595 & 0.475 & 1.945 \\
  P-n55-k8-260 & 0.031 & 0.078 & 0.129 & 0.324 & 0.496 & 0.780 \\
  P-n60-k10-24 & 0.035 & 0.080 & 0.114 & 0.272 & 0.480 & 0.901 \\
  P-n60-k15-8 & 0.014 & 0.017 & 0.024 & 0.028 & 0.034 & 0.036 \\
  P-n65-k10-102 & 0.023 & 0.063 & 0.058 & 0.196 & 0.140 & 0.548 \\
  P-n70-k10-12 & 0.188 & 0.396 & 0.654 & 2.159 & 2.887 & 11.198 \\
  P-n76-k4-41 & 46.847 & 24.615 & \textit{t/o} & \textit{t/o} & \textit{t/o} & \textit{t/o} \\
  P-n76-k5-16 & 10.934 & 8.169 & 55.936 & 94.270 & \textit{t/o} & \textit{t/o} \\
  \midrule
  \multicolumn{7}{l}{\textbf{roberti}} \\
  E-n101-k14\_a & 0.148 & 0.441 & 0.773 & 1.532 & 5.413 & 5.288 \\
  E-n101-k14\_b & 0.121 & 0.365 & 0.659 & 1.290 & 4.616 & 4.297 \\
  E-n101-k8\_a & 0.361 & 1.648 & 2.410 & 11.429 & 38.492 & 117.529 \\
  E-n101-k8\_b & 0.154 & 0.752 & 0.944 & 3.008 & 7.878 & 10.266 \\
  E-n76-k10\_a & 0.088 & 0.256 & 0.519 & 1.281 & 2.913 & 4.350 \\
  E-n76-k10\_b & 0.054 & 0.194 & 0.256 & 0.819 & 1.025 & 2.402 \\
  E-n76-k14\_a & 0.039 & 0.092 & 0.153 & 0.265 & 0.617 & 0.688 \\
  E-n76-k14\_b & 0.040 & 0.092 & 0.159 & 0.247 & 0.651 & 0.714 \\
  E-n76-k7\_a & 0.140 & 0.640 & 1.308 & 4.596 & 13.841 & 32.507 \\
  E-n76-k7\_b & 0.117 & 0.478 & 0.622 & 3.497 & 3.637 & 17.070 \\
  E-n76-k8\_a & 0.099 & 0.416 & 0.629 & 2.517 & 4.124 & 12.644 \\
  E-n76-k8\_b & 0.054 & 0.211 & 0.158 & 0.899 & 0.546 & 2.871 \\
  F-n135-k7\_a & 3.620 & 64.392 & \textit{t/o} & \textit{t/o} & \textit{t/o} & \textit{t/o} \\
  F-n45-k4\_a & 0.125 & 0.317 & 0.415 & 2.829 & 3.646 & 51.195 \\
  F-n72-k4\_a & 9.742 & 3.115 & 15.717 & 74.461 & 49.041 & \textit{t/o} \\
  M-n121-k7\_a & 2.903 & \textit{t/o} & \textit{t/o} & \textit{t/o} & \textit{t/o} & \textit{t/o} \\
  M-n121-k7\_b & 1.937 & 52.127 & \textit{t/o} & \textit{t/o} & \textit{t/o} & \textit{t/o} \\
  M-n151-k12\_a & 0.674 & 4.964 & 5.929 & 32.903 & 102.890 & \textit{t/o} \\
  M-n151-k12\_b & 0.591 & 4.140 & 5.489 & 25.701 & 86.183 & \textit{t/o} \\
  M-n200-k16\_a & 1.612 & 14.223 & 16.548 & 83.345 & \textit{t/o} & \textit{t/o} \\
  M-n200-k16\_b & 1.454 & 11.413 & 15.165 & 61.544 & \textit{t/o} & \textit{t/o} \\
  M-n200-k17\_a & 1.332 & 10.727 & 11.460 & 58.021 & \textit{t/o} & \textit{t/o} \\
  M-n200-k17\_b & 1.332 & 9.988 & 12.033 & 51.121 & \textit{t/o} & \textit{t/o} \\
  P-n101-k4\_a & 0.561 & 4.468 & 6.561 & 53.694 & \textit{t/o} & \textit{t/o} \\
  P-n101-k4\_b & 0.266 & 1.317 & 1.473 & 7.599 & 18.415 & 58.645 \\
  P-n70-k10\_a & 0.038 & 0.115 & 0.136 & 0.379 & 0.359 & 1.236 \\
  P-n70-k10\_b & 0.053 & 0.103 & 0.091 & 0.335 & 0.285 & 1.172 \\
  P-n76-k4\_a & 0.225 & 1.202 & 2.555 & 14.021 & 38.199 & \textit{t/o} \\
  P-n76-k4\_b & 0.103 & 0.463 & 0.678 & 3.576 & 6.346 & 27.046 \\
  P-n76-k5\_a & 0.161 & 0.902 & 1.724 & 9.157 & 21.948 & 117.527 \\
  P-n76-k5\_b & 0.097 & 0.422 & 0.527 & 2.655 & 3.675 & 15.134 \\
\end{longtable}
\end{footnotesize}

\begin{footnotesize}
\begin{longtable}{lrrrrrr}
  \caption{Per-instance wall-clock time (s) on the 56 Solomon rcspp instances: \texttt{bgspprc} (\texttt{para\_bidir}, \emph{bg}) against parallel pull labelling (\emph{pull}). \textit{t/o} marks a 120\,s timeout. Per-column shifted geometric means reproduce \autoref{tab:rcspp}.}\label{tab:per-instance-pull}\\
  \toprule
  Instance & \multicolumn{2}{c}{$n_g=8$} & \multicolumn{2}{c}{$n_g=16$} & \multicolumn{2}{c}{$n_g=24$} \\
   & bg & pull & bg & pull & bg & pull \\
  \midrule
  \endfirsthead
  \multicolumn{7}{l}{\small\itshape (\autoref{tab:per-instance-pull} continued)}\\
  \toprule
  Instance & \multicolumn{2}{c}{$n_g=8$} & \multicolumn{2}{c}{$n_g=16$} & \multicolumn{2}{c}{$n_g=24$} \\
   & bg & pull & bg & pull & bg & pull \\
  \midrule
  \endhead
  \midrule \multicolumn{7}{r}{\small\itshape continued on next page}\\
  \endfoot
  \bottomrule
  \endlastfoot
  C101 & 0.029 & 0.010 & 0.029 & 0.010 & 0.030 & 0.010 \\
  C102 & 0.205 & 0.017 & 0.218 & 0.016 & 0.211 & 0.017 \\
  C103 & 0.423 & 0.023 & 0.431 & 0.026 & 0.452 & 0.030 \\
  C104 & 0.761 & 0.057 & 0.894 & 0.319 & 0.812 & 0.121 \\
  C105 & 0.060 & 0.009 & 0.067 & 0.009 & 0.070 & 0.009 \\
  C106 & 0.099 & 0.010 & 0.083 & 0.010 & 0.082 & 0.010 \\
  C107 & 0.113 & 0.009 & 0.125 & 0.009 & 0.123 & 0.010 \\
  C108 & 0.146 & 0.013 & 0.139 & 0.012 & 0.138 & 0.013 \\
  C109 & 0.254 & 0.017 & 0.240 & 0.017 & 0.256 & 0.017 \\
  C201 & 0.242 & 0.010 & 0.282 & 0.010 & 0.269 & 0.010 \\
  C202 & 6.180 & 0.030 & 6.264 & 0.040 & 6.884 & 0.084 \\
  C203 & 17.696 & 0.163 & 20.527 & 2.478 & 18.696 & 0.571 \\
  C204 & 32.341 & 2.961 & 38.406 & 16.294 & 94.144 & \textit{t/o} \\
  C205 & 0.613 & 0.020 & 0.600 & 0.020 & 0.593 & 0.017 \\
  C206 & 1.077 & 0.037 & 1.057 & 0.033 & 1.105 & 0.030 \\
  C207 & 1.736 & 0.211 & 1.729 & 0.277 & 1.634 & 0.095 \\
  C208 & 1.510 & 0.044 & 1.516 & 0.040 & 1.557 & 0.041 \\
  R101 & 0.015 & 0.006 & 0.014 & 0.006 & 0.013 & 0.006 \\
  R102 & 0.107 & 0.011 & 0.108 & 0.012 & 0.084 & 0.012 \\
  R103 & 0.201 & 0.019 & 0.207 & 0.019 & 0.217 & 0.018 \\
  R104 & 0.514 & 0.038 & 0.662 & 0.063 & 0.664 & 0.056 \\
  R105 & 0.046 & 0.008 & 0.042 & 0.008 & 0.053 & 0.008 \\
  R106 & 0.171 & 0.017 & 0.150 & 0.017 & 0.166 & 0.018 \\
  R107 & 0.346 & 0.037 & 0.354 & 0.042 & 0.358 & 0.038 \\
  R108 & 0.656 & 0.063 & 0.715 & 0.073 & 0.778 & 0.081 \\
  R109 & 0.117 & 0.013 & 0.137 & 0.013 & 0.120 & 0.013 \\
  R110 & 0.252 & 0.020 & 0.267 & 0.020 & 0.255 & 0.020 \\
  R111 & 0.310 & 0.037 & 0.338 & 0.041 & 0.342 & 0.042 \\
  R112 & 0.581 & 0.047 & 0.610 & 0.062 & 0.600 & 0.063 \\
  R201 & 0.738 & 0.016 & 0.769 & 0.017 & 0.724 & 0.017 \\
  R202 & 5.760 & 0.084 & 5.695 & 0.101 & 5.931 & 0.139 \\
  R203 & 18.320 & 0.293 & 19.519 & 0.816 & 34.945 & 7.891 \\
  R204 & 36.713 & 1.394 & 57.020 & 5.782 & \textit{t/o} & \textit{t/o} \\
  R205 & 2.269 & 0.070 & 2.442 & 0.105 & 2.375 & 0.092 \\
  R206 & 8.736 & 0.184 & 10.946 & 0.589 & 13.767 & 1.706 \\
  R207 & 21.125 & 0.754 & 29.556 & 2.372 & 64.479 & 12.027 \\
  R208 & 38.875 & 2.124 & 62.849 & 30.547 & \textit{t/o} & \textit{t/o} \\
  R209 & 5.858 & 0.130 & 6.923 & 0.274 & 7.601 & 0.296 \\
  R210 & 7.074 & 0.178 & 8.841 & 0.515 & 12.386 & 1.392 \\
  R211 & 12.821 & 0.452 & 23.624 & 1.852 & 72.988 & 11.193 \\
  RC101 & 0.035 & 0.008 & 0.032 & 0.008 & 0.027 & 0.008 \\
  RC102 & 0.097 & 0.014 & 0.088 & 0.014 & 0.110 & 0.016 \\
  RC103 & 0.240 & 0.024 & 0.260 & 0.026 & 0.253 & 0.025 \\
  RC104 & 0.628 & 0.075 & 0.813 & 0.262 & 0.800 & 0.230 \\
  RC105 & 0.067 & 0.012 & 0.069 & 0.014 & 0.071 & 0.014 \\
  RC106 & 0.085 & 0.013 & 0.081 & 0.013 & 0.100 & 0.013 \\
  RC107 & 0.248 & 0.022 & 0.243 & 0.023 & 0.238 & 0.022 \\
  RC108 & 0.442 & 0.059 & 0.533 & 0.059 & 0.526 & 0.060 \\
  RC201 & 0.657 & 0.018 & 0.672 & 0.019 & 0.652 & 0.018 \\
  RC202 & 4.291 & 0.062 & 4.426 & 0.061 & 4.232 & 0.059 \\
  RC203 & 16.025 & 0.744 & 18.433 & 0.924 & 18.716 & 1.158 \\
  RC204 & 44.404 & 10.969 & 60.759 & 88.971 & \textit{t/o} & 77.071 \\
  RC205 & 1.934 & 0.038 & 2.067 & 0.047 & 2.096 & 0.038 \\
  RC206 & 1.968 & 0.058 & 2.172 & 0.063 & 2.182 & 0.071 \\
  RC207 & 4.864 & 0.162 & 6.617 & 0.515 & 8.367 & 0.784 \\
  RC208 & 14.792 & 0.826 & 90.205 & 14.600 & \textit{t/o} & 35.112 \\
\end{longtable}
\end{footnotesize}

\bibliographystyle{unsrt}
\bibliography{references}

\end{document}